\documentclass[11pt, a4paper]{article}
\usepackage{latexsym}
\usepackage{indentfirst}
\usepackage{graphicx}
\usepackage{amsmath}
\usepackage{amssymb}
\begin{document}
\title{Stability of equilibrium and periodic solutions of a
delay equation modeling leukemia}
\author{Anca-Veronica Ion$^1$, Raluca-Mihaela Georgescu$^2$\\
$^1$"Gh. Mihoc-C. Iacob" Institute of Mathematical Statistics\\
and Applied Mathematics of the Romanian Academy, Bucharest, Romania\\
$^2$University of Pite\c sti, Romania}

\date{}
\maketitle

\begin{abstract}
We consider a delay differential equation that occurs in the study
of chronic myelogenous leukemia. The equation was investigated
numerically in \cite{PM-M}, \cite{PM-B-M}. In \cite{I09} we studied
in detail  the stability of the two equilibrium points.  In the
present work, we shortly remind the results of \cite{I09} and then
concentrate on the study of stability of periodic solutions emerged
by Hopf bifurcation from
 the non-trivial equilibrium point. We give the algorithm for
approximating a center manifold at a typical point of Hopf
bifurcation (and  an unstable manifold in the vicinity of such a
point, where such a manifold exists). Then we
 find the normal form of the equation restricted to the
center manifold, by computing the first Lyapunov coefficient. The
normal form allows us to establish the stability properties of the
periodic solutions occurred by Hopf bifurcation.

\textbf{Acknowledgement.} Work supported by Grant 12/11.02.2008
within the framework of the Russian Foundation for
Basic Research - Romanian Academy collaboration.
\end{abstract}

\section{The problem}
We consider the delay differential equation that occurs in the study
of periodic chronic myelogenous leukemia \cite{PM-M}, \cite{PM-B-M}
\begin{equation}\label{eq}
\dot{x}(t)=-\left[\frac{\beta_0}{1+x(t)^n}+\delta\right]x(t)+k\frac{\beta_0x(t-r)}{1+x(t-r)^n}.
\end{equation}
Here  $\beta_0,\,n,\,\delta,\,\gamma,\,r$ are positive parameters
and $k=2e^{-\gamma r}$.  The unknown function, $x(\cdot),$ should
also be nonnegative, being a non-dimensional density of cells. We do
not insist on the significance of the function $x(.)$ or in that of
the parameters, since these are extensively presented in
\cite{PM-M}, \cite{PM-B-M}. Anyway, inasmuch as  the equation
represents a good model for the periodic chronic myelogenous
leukemia, both the equilibrium solutions and the periodic solutions
are important, as are their stability properties.

 In the following we use the space
 $\mathcal{B}=C([-r,0],\mathbb{R})$ (the space of continuous, real valued functions
 defined on $[-r,0],$ with the supremum norm, denoted by $|x|_0$).
Given a function $x:[-r,T)\mapsto \mathbb{R},\,T>0$ and a $0\leq t <
T,$ we define the function $x_t\in \mathcal{B}$ by $x_t(s)=x(t+s).$

Equation (\ref{eq}) may be written as \begin{equation}\label{abseq}
\dot{x}=h(x_t),
\end{equation}
where $h:\mathcal{B}\mapsto\mathbb{R},$ and  we impose to this
equation the initial condition
\begin{equation}\label{ci}x_0=\phi\in \mathcal{B}.
\end{equation}

In \cite{I09} we proved that problem (\ref{abseq}), (\ref{ci}) has
an unique, defined on $[-r,\infty)$  bounded solution.

\subsection{Equilibrium solutions \cite{PM-M}, \cite{PM-B-M}}

As shown in \cite{PM-M}, the equilibrium points of the problem are
$$x_1=0,\,\, x_2=(\frac{\beta_0}{\delta}(k-1)-1)^{1/n}.$$ The second
one is acceptable from the biological point of view if and only if
is strictly positive that is, if and only if
\begin{equation}\label{cond} \frac{\beta_0}{\delta}(k-1)-1
> 0.\end{equation}

In terms of $r,$ the above inequality may be written as
\[r<r_{max}:=-\frac{1}{\gamma}\ln\frac{1}{2}\left(1+\frac{\delta}{\beta_0}\right),
\]
and since the delay $r$ is positive, the condition
$\delta/\beta_0<1$ follows.

The biological interpretation of function
$\beta(x)=\beta_{0}/(1+x^n)$ \cite{PM-M} shows that the condition
$\beta(x_2)=\delta/(k-1)>0$ should be fulfilled. This is equivalent
to $k>1.$

The linearized equation around one of the equilibrium points is
\begin{equation}\label{lineq}
\dot{z}(t)=-[B_1+\delta]z(t)+kB_1z(t-r),\end{equation} with
$B_1=\beta'(x^*)x^*+\beta(x^*),\,x^*=x_1$ or $x^*=x_2.$

The eigenvalues of the infinitesimal generator of the semigroup of
operators generated by equation (\ref{lineq}) are the solutions of
the characteristic equation
\begin{equation}\label{chareq}
\lambda+\delta+B_1=kB_1e^{-\lambda r}.
\end{equation}

\section{Stability of equilibrium solution}

Some conclusions on the stability of equilibrium points are
presented in \cite{PM-M}, \cite{PM-B-M}. The study in \cite{I09} has
lead us to the following conclusions (obtained by using the results
in \cite{Ha}), conclusions that bring some corrections to the
stability results in \cite{PM-M}, \cite{PM-B-M}.

\subsection{Stability of $x_1$}
In this case $B_1=\beta_0.$ The equilibrium point  $x_1$ is stable
as long as it is the single equilibrium point, that is when
$\frac{\beta_0}{\delta}(k-1)\leq 1$. When the second equilibrium
point occurs ($\frac{\beta_0}{\delta}(k-1)> 1$), $x_1$ becomes
unstable. For $\frac{\beta_0}{\delta}(k-1)=1,$ equation
(\ref{chareq})  admits the solution $\lambda=0.$ Hence the change of
stability occurs by traversing the eigenvalue $\lambda=0.$

\subsection{Stability of  $x_2$}

We still rely on \cite{I09} in this subsection.

For $x_2$,
\begin{equation}\label{B1}
B_1=\beta_0[n-(n-1)A]/A^2\end{equation} where
$A=\beta_0(k-1)/\delta.$

Let us define for future use
\[r_n=-\frac{1}{\gamma}\ln \left\{\frac{1}{2}\left(\frac{\delta}{\beta_0}\frac{n}{n-1}+
1\right)\right\}.
\]

We have the following distinct situations.

\textbf{I.} $B_1<0$ that is equivalent to
$\frac{\beta_0}{\delta}(k-1)>\frac{n}{n-1},$ and, in terms of $r,$
equivalent to $0\leq r\leq r_n.$ Two subcases are to be considered.

\textbf{I.A.} $B_1<0$ and $\delta+B_1<0.$ In this case $Re \lambda
<0$ for all eigenvalues $\lambda$ if and only if
$|\delta+B_1|<|kB_1|$ and
\begin{equation}\label{A}
\frac{\arccos{((\delta+B_1)/kB_1)}}{\omega_0}<r<\frac{1}{|\delta+B_1|},
\end{equation}
where $\omega_0$ is the solution in $(0,\pi/r)$ of the equation
$\omega\cot (\omega r)=-(\delta+B_1).$ Alternately, $\omega_0$ may
be defined as
\begin{equation}\label{om0}\omega_0=\frac{1}{r}T^{-1}(-(\delta+B_1)r).\end{equation}
where $T:[0,\pi)\mapsto (-\infty,\,1]$ is given by
\begin{equation}\label{T}
T(y)=\left\{%
\begin{array}{ll}
    y\cot(y), &  y\in (0,\pi); \\
    1, & y=0. \\
\end{array}%
\right.
\end{equation}

Hence, when the above conditions are satisfied, the equilibrium
solution $x_2$ is stable.

When $r|\delta+B_1|=1,$ the solution $x_2$ is unstable \cite{I09}.

\textbf{I.B.} $B_1<0$ and $\delta+B_1>0.$

In this case, $Re\lambda <0$ for all eigenvalues $\lambda$ if and
only if
\begin{equation}\label{B}
\delta+B_1>|kB_1|\,\,\mathrm{or}\,\left\{\delta+B_1\leq |kB_1|\,\,
\mathrm{and}\,\,r<\frac{\arccos{((\delta+B_1)/kB_1)}}{\omega_0}\right\}
\end{equation}
where $\omega_0$ is defined as above. When these conditions are
satisfied, $x_2$ is stable.

\textbf{Remarks. a)} If $B_1<0,\,\delta+B_1=0,$ the solution is
stable if and only if $-kB_1 r<\pi/2$ while, for this case, the
point $kB_1 r =-\pi/2$ is a Hopf bifurcation point.

\textbf{b)} Assume that we vary $r$ and keep all the other
parameters fixed. The conditions (\ref{A}) or (\ref{B}) for $r$ are
not as simple as they seem, because $B_1$ is itself a function of
$r,$ (being a function of $k$). Let us consider the function
\begin{equation}\label{h}
g(r)=T^{-1}(-(\delta+B_1(r))r)-\arccos\left(\frac{\delta+B_1(r)}{k(r)B_1(r)}\right).
\end{equation}

If for a certain $r^*$ we have $g(r^*)=0$ (that is the condition for
the change of stability), in order to find whether a value $r_1$ in
a neighborhood of $r^*$
 is in the stability zone or not, we have
to know the sign of $g(r_1)$, hence we have to study the monotony
properties of function $g$ in a neighborhood of $r^*$. This will be
done in the last section for an example.

\textbf{II.} $B_1>0$ that is equivalent to
$\frac{\beta_0}{\delta}(k-1)<\frac{n}{n-1},$ and, in terms of $r,$
equivalent to $r_n\leq r\leq r_{max}.$

In this case, we can only have $\delta+B_1>0,$ and, by \cite{Ha},
$Re\lambda<0$ for all eigenvalues $\lambda$ if and only if
$kB_1<\delta+B_1.$

By using (\ref{B1}) we see the above inequality is equivalent to
$\frac{\beta_0(k-1)}{\delta}>1$ that is already satisfied, since
$x_2$ exists and is positive. It follows that if $B_1>0$ then $x_2$
is stable.

\section{Hopf bifurcation points}

By denoting $\lambda=\mu+i\omega,$ and by equating the real,
respectively the imaginary part of the characteristic equation
(\ref{chareq}), we obtain
\begin{equation}\label{re,im}
\mu+\delta+B_1=kB_1 e^{-\mu r}\cos(\omega r),
\end{equation}
\[\omega=-kB_1 e^{-\mu r}\sin(\omega r).
\]
It is useful to consider the case $\mu=0$ in the above equations,
\begin{equation}\label{muegal0}
\delta+B_1=kB_1 \cos(\omega r),
\end{equation}
\[\omega=-kB_1 \sin(\omega r).
\]
In the above stability discussion, in \textbf{I.}, the case
\begin{equation}\label{r-H}
r=\frac{\arccos((\delta+B_1)/(kB_1))}{\omega_0}
\end{equation}
occurs on the frontier of the stability domain. This relation,
together with $\omega_0\cot(\omega_0 r)=-(\delta+B_1)$ imply that
$\omega_0=\sqrt{(kB_1)^2-(\delta+B_1)^2}$ and that the pair
$\mu^*=0,\,\omega^*=\omega_0$ represents a solution of
(\ref{re,im}).

Hence a point in the parameter space, where relation (\ref{r-H}) is
satisfied, might be a Hopf bifurcation point when one of the
parameters is varied. The Hopf bifurcation may be either a
non-degenerated  or a degenerated one.

To be more precise, we consider the vector of parameters
$(\beta_0,\,n,\,\delta,\,\gamma,\,r)$ and denote it by
$\mathbf{a}=(a_i)_{1\leq i\leq 5}.$

We choose a point in the parameter space,
$\mathbf{a}^*=(a^*_i)_{1\leq i\leq 5},$ such that, for this choice
of parameters, there are two pure imaginary eigenvalues
$\widetilde{\lambda}_{1,2}(\mathbf{a}^*)=\pm i
\omega^*,\,\omega^*:=\omega(\mathbf{a}^*)>0$   and all other
eigenvalues $\lambda$ have strictly negative real parts (i.e. a
point where condition (\ref{r-H}) holds).

In the following we keep fixed all parameters excepting one, $a_j,$
that we vary\,($j$ fixed, $1\leq j\leq 5,$). For further simplicity,
we denote $a_j:=\alpha,$ such that
$\widetilde{\lambda}_{1,2}(\mathbf{a})$ becomes
$\widetilde{\lambda}_{1,2}(\alpha)$ and such that $Re
\widetilde{\lambda}_{1,2}(\alpha^*)=0$.

If $\frac{\partial Re \lambda_i}{\partial \alpha}(\alpha^{*})\neq
0,\,i=1,2,$ then at $\alpha^*$ a Hopf bifurcation occurs and, in a
neighborhood $U^*$ of $\alpha^*,$ either to the left or to the right
of $\alpha^*$ we have $Re \widetilde{\lambda}_{1,2}(\alpha)<0,$
while $Re \widetilde{\lambda}_{1,2}(\alpha)>0$ at the opposite side
of $\alpha^*.$ The neighborhood $U^*$ is such that for all
$\alpha\in U^*$ all other eigenvalues besides
$\widetilde{\lambda}_{1,2}$ have real parts less than a fixed
negative value.

\section{On the stability properties of the periodic solutions emerged by Hopf bifurcation}

We consider a periodic solution occured by Hopf bifurcation, as
described in the previous section. In order to establish whether or
not such a periodic solution is stable, we have to compute the
normal form of the reduced to the center manifold (or to the
unstable manifold) problem. In this section we present the algorithm
for obtaining this normal form and the consequences on the stability
of periodic solutions.

Since we study local bifurcations around $x_2,$ throughout this
section we work with eq. (\ref{eq}) translated by $x_2,$ that is,
the new unknown function $z=x-x_2$ is considered.

\subsection{The delay equation as an equation in a Banach space}

We assume that we are in one of the situations, mentioned in Section
3, where a Hopf bifurcation occurs. In order to construct the center
(or the unstable) manifold at the equilibrium point $x_2$  we need
to write our equation as a differential equation in a Banach space.
We use the method of Teresa Faria, from \cite{Far}, as we did also
in \cite{I04}, \cite{I07}. We consider the space
\[\mathcal{B}_0=\left\{\psi:[-r,0]\mapsto \mathbb{R},\, \psi\, \mathrm{is\, continuous\, on\,}[-r,0)\wedge\,\exists
\lim_{s\rightarrow 0}\psi(s)\in \mathbb{R}  \right\},
\]
that is formed of functions $\psi=\varphi+\sigma d_{0},$ where
$\varphi\in \mathcal{B},$ $\sigma\in \mathbb{R}$ and\\
$d_0:[-r,0]\mapsto \mathbb{R},$
$$d_0(s)=\left\{\begin{array}{cc}
                0, & s\in[-r,0), \\
                1, & s=0. \\
              \end{array}\right.
$$
The space $\mathcal{B}_0$ is normed by
$\|\psi\|=|\varphi|_0+|\sigma|.$ The linearized equation
(\ref{lineq}) can be written as
\begin{equation}\label{lineq1}
\dot{z}=Lz_{t},\end{equation}
\[\quad L\varphi =\int_{-r}^0d\eta(\theta)\varphi(\theta),
\]
(the integral is a Stieljes one) with
\[\eta(s)=\left\{\begin{array}{c}
              -kB_1,\quad s=-r; \\
              0,\quad s\in (-r,0); \\
              -(B_1+\delta),\quad s=0.\\
            \end{array}\right.
\]
Nonlinear equation (\ref{eq}) may be written as
\begin{equation}\label{nonlineq1}\dot{z}(t)=L(z_t)+f(z(t),z(t-r)),
\end{equation}
where $f(\cdot,\cdot)$ comprises the nonlinear terms of the
equation. Since $z(t)$ and $z(t-r)$ are functions of $z_t,$ eq.
(\ref{nonlineq1}) may be written as
\begin{equation}\label{nonlineq2}\dot{z}(t)=L(z_t)+\widetilde{f}(z_t).
\end{equation}
In \cite{Far}, the linear operator
$A:C^{1}([-r,0],\mathbb{R})\subset\mathcal{B}_0\mapsto
\mathcal{B}_0$
\begin{equation}A(\varphi)=\dot{\varphi}+d_0[L(\varphi)-\dot{\varphi}(0)]
\end{equation}
is defined and it is proved that this is the infinitesimal generator
of the semigroup of operators $\{S(t)\}_{t\geq 0}$ given by
$S(t)(\phi)=z_t(\phi),$ where $z(t,\phi)$ is the solution of
equation (\ref{lineq1}) with the initial condition $z_0=\phi.$ Then
the nonlinear equation may be written as an equation in
$\mathcal{B}_0,$ that is
\begin{equation}\label{nonlineq3}\frac{dz_t}{dt}=A(z_t)+d_0\widetilde{f}(z_t).
\end{equation}
\subsection{Space of eigenfunctions corresponding to $\widetilde{\lambda}_{1,2}$,
the projector on this space} We return for the moment to the linear
equation (\ref{lineq}). The eigenfunctions corresponding to the two
 eigenvalues
$\widetilde{\lambda}_{1,2}=\mu(\alpha) \pm i \omega(\alpha)$
 are given by
$\varphi_{1,2}(s)=e^{\widetilde{\lambda}_{1,2} s},\,s\in [-r,0].$
Since the eigenfunctions are complex functions, we need to use the
complexification of the spaces $\mathcal{B},\,\mathcal{B}_0$ that we
denote by $\mathcal{B}_{C},\,\mathcal{B}_{0C},$ respectively. We
denote by $\mathcal{M}$ the subspace of $\mathcal{B}_C$ generated by
$\varphi_{1,2}(\cdot).$

In \cite{Far}, for an equation of the form \eqref{nonlineq3}, a
projector defined on $\mathcal{M}$ is constructed. For its
construction we need a bilinear form that is defined by using the
transposed equation of the linearized equation. The transposed
equation is
\[\dot{y}(s)=-\int_{-r}^0y(s-\theta)d[\eta(\theta)],
\]
hence
\[\dot{y}(s)=(B_1+\delta)y(s)-kB_1y(s+r).
\]

The characteristic equation for the adjoint problem is
\[\lambda=B_1+\delta-kB_1e^{\lambda r},\] that admits the solutions $\lambda_{1,2}=-\mu(\alpha) \pm
i\omega(\alpha).$ The corresponding eigenfunctions are
$\psi_1(\zeta)=e^{(-\mu+i\omega)\zeta},\,\psi_2(\zeta)=e^{(-\mu-i\omega)\zeta},$
$\zeta\in [0,r]$.

In order to construct the projector defined on $\mathcal{B}_{0C}$
with values on $\mathcal{M},$ we define the bilinear form \cite{HaL}
\[\langle\psi,\varphi\rangle=\psi(0)\varphi(0)-\int_{-r}^0\int_{0}^\theta
\psi(\zeta-\theta)d\eta(\theta)\varphi(\zeta)d\zeta=
\]
\[=\psi(0)\varphi(0)+
kB_1\int_{-r}^0\psi(\zeta+r)\varphi(\zeta)d\zeta,
\]
where $\psi:[0,r]\rightarrow \mathbb{R}.$

We look for linear combinations of the vectors $\psi_j,$ that we
will denote $\Psi_i,\,i=1,2,$ such that
$\langle\Psi_i,\,\varphi_j\rangle=\delta_{ij}.$ For this we
determine the $2\times 2$ matrix $E,$ with elements
$e_{ij}=\langle\psi_i,\,\varphi_j\rangle$.
\[e_{11}=(\psi_1,\varphi_1)=\psi_1(0)\varphi_1(0)+
kB_1\int_{-r}^0e^{(-\mu+i\omega)(\theta+r)}e^{(\mu+i\omega)\theta}d\theta=1-e^{\mu
r},
\]
\[e_{12}=(\psi_1,\varphi_2)=\psi_1(0)\varphi_2(0)+
kB_1\int_{-r}^0e^{(-\mu+i\omega)(\theta+r)}e^{(\mu-i\omega)\theta}d\theta=1+(\delta+B_1+\mu-i\omega)r,
\]
\[e_{21}=(\psi_2,\varphi_1)=\psi_2(0)\varphi_1(0)+
kB_1\int_{-r}^0e^{(-\mu-i\omega)(\theta+r)}e^{(\mu+i\omega)\theta}d\theta=1+(\delta+B_1+\mu+i\omega)r,
\]
\[e_{22}=(\psi_2,\varphi_2)=\psi_2(0)\varphi_2(0)+
kB_1\int_{-r}^0e^{(-\mu-i\omega)(\theta+r)}e^{(\mu-i\omega)\theta}d\theta=1-e^{\mu
r}.
\]
We have $\det E = e_{11}e_{22}-e_{12}e_{21}=(1-e^{\mu
r})^2-[(1+(\delta+B_1+\mu)r)^2+\omega^2r^2]. $ Since
\[\left(%
\begin{array}{c}
  \Psi_1 \\
  \Psi_2 \\
\end{array}%
\right)=E^{-1}\left(%
\begin{array}{c}
  \psi_1 \\
  \psi_2 \\
\end{array}%
\right),
\]
we obtain
\[\Psi_1=\frac{1}{\det E}[e_{22}\psi_1-e_{12}\psi_2 ]=\frac{-e_{12}}{\det
E}\psi_2,\quad \Psi_2=\overline{\Psi}_1.
\]

We can now write the projector defined on $\mathcal{B}_{0C}$ and
with values in $\mathcal{M}.$ It is given, for
$\psi=\phi+d_0\sigma\in\mathcal{B}_{0C}$ by
\begin{equation}\mathcal{P}(\psi)=\left(\langle\Psi_1,\phi\rangle+\Psi_1(0)\sigma\right)\varphi_1
+\left(\langle\Psi_2,\phi\rangle+\Psi_2(0)\sigma\right)\varphi_2.\end{equation}
If $\sigma=0$ and, thus, $\psi=\phi\in\mathcal{B}_C,$ we have
$\mathcal{P}(\phi)=\langle\Psi_1,\phi\rangle\varphi_1
+\langle\Psi_2,\phi\rangle\varphi_2.$

Now, for the solution $z(\cdot,\phi)$ of (\ref{nonlineq3}) we can
write
\[z_t=\phi_1u_1(t)+\phi_2u_2(t)+\mathbf{v}(t),\,
\]
with $u_1(t)=\langle\Psi_1,\,z_t
\rangle,\,\,u_2(t)=\langle\Psi_2,\,z_t
\rangle,\,\,\mathbf{v}=(I-\mathcal{P})z_t.$

We project now the equation (\ref{nonlineq3}) by $\mathcal{P},$ to
obtain,
\begin{equation}\label{proj-eq}
\frac{d\mathbf{u}}{dt}=B\mathbf{u}+\Psi(0)\widetilde{f}(\phi_1u_1+\phi_2u_2+\mathbf{v}),
\end{equation}
where $$\mathbf{u}=\left(%
\begin{array}{c}
  u_1 \\
  u_2 \\
\end{array}%
\right), \,\,B=\left(%
\begin{array}{cc}
  \widetilde{\lambda}_1 & 0 \\
  0 & \widetilde{\lambda}_2 \\
\end{array}\right),\,\,\Psi=\left(%
\begin{array}{c}
  \Psi_1 \\
  \Psi_2 \\
\end{array}%
\right).$$

Remark that $u_1,u_2$ are complex functions, but, since our initial
problem is formulated in terms of real functions, we have
$u_1=\overline{u}_2,$ and the two scalar equations comprised in
(\ref{proj-eq}) are complex conjugated one to the other. Hence it
suffices to study one of the them, let us say the equation for
$u_1.$ We denote $u_1=u$ and the equation to be studied in the
following is
\begin{equation}\label{proj-eq1}
\frac{du}{dt}=(\mu+i\omega
)u+\Psi_1(0)\widetilde{f}(\varphi_1u+\varphi_2\overline{u}+\mathbf{v}).
\end{equation}
For $\alpha=\alpha^*$ we have $\mu(\alpha^*)=0$ hence
\begin{equation}\label{psi}\Psi_1(0)=\frac{1+(\delta+B_1-i\omega)r}{[1+(\delta+B_1)r]^2+\omega^2r^2} .
\end{equation}

\subsection{The equation restricted to the invariant (center or unstable) manifold}

For those $\alpha \in U^*$ where $Re
\widetilde{\lambda}_{1,2}(\alpha) \geq 0,$ since all other
eigenvalues have negative real part, there is a local invariant
manifold, the local unstable manifold, $W_{loc}^{u},$ for $Re
\widetilde{\lambda}_{1,2}(\alpha) > 0,$ or the local center
manifold, $W_{loc}^{c},$ at $\alpha^*$ ($Re
\widetilde{\lambda}_{1,2}(\alpha^*) = 0$).

The local invariant manifold is an at least $C^1$  invariant
manifold, tangent to the space $\mathcal{M}$ at the point $z=0$
(that is $x=x_2$), and it is the graph of a function
$w(\alpha)(\cdot)$ defined on a neighborhood of zero in
$\mathcal{M}$ and taking values in $(I-\mathcal{P})\mathcal{B_{C}.}$
A point on the local invariant manifold has the form
$u\,\varphi_1+\overline{u}\,\varphi_2+w(\alpha)(u\,\varphi_1+\overline{u}\,\varphi_2).$
Since it is an invariant manifold, if the initial condition is taken
on the manifold, then its image through the semigroup
$\{\widetilde{T}(t)\}_{t\geq 0},$
$\widetilde{T}(t)(\phi)=z_t(\phi)$, is still on the manifold
($\widetilde{T}(t)(\phi)={T}(t)(\phi)-x_2,$ where $T(t)$ is defined
at the end of Subsection 1.1). Hence
\begin{equation}\label{zt}\widetilde{T}(t)(\phi)=u(t)\varphi_1+\overline{u(t)}\varphi_2+
w(\alpha)(u(t)\varphi_1+\overline{u(t)}\varphi_2),
\end{equation}
with $u(\cdot),$ solution of the equation

\begin{equation}\label{proj-eq1}
\frac{du}{dt}=\widetilde{\lambda}_{1}(\alpha)u+\Psi_1(0)\widetilde{f}(\varphi_1u+\varphi_2\overline{u}+w(\alpha)(u\varphi_1+
\overline{u}\varphi_2)),
\end{equation}
with the initial condition $u(0)=u_0,$ where
$\mathcal{P}(\phi)=u_0\varphi_1+\overline{u}_0\varphi_2.$ The real
and the imaginary parts of this complex equation, represent the
two-dimensional restricted to the center manifold problem. We can
study this problem with the tools of planar dynamical systems theory
(see, e.g. \cite{K}).

\subsection{The normal form of the reduced equation}

In order to find the normal form of equation (\ref{proj-eq1}) we
must consider the series of powers of $u,\,\overline{u}$ for
$w(\alpha)$ and for $\widetilde{f}(z_t).$ We set
\begin{equation}\label{seriesw}\widetilde{w}(\alpha)(u,\overline{u})=
\sum_{j+k\geq2}\frac{1}{j!k!}w_{jk}(\alpha)u^{j}\overline{u}^{k},
\end{equation}
where
$\widetilde{w}(\alpha)(u,\overline{u}):=w(\alpha)(u(t)\varphi_1+\overline{u(t)}\varphi_2),$
 $w_{jk}(\alpha)\in\mathcal{B},$  $u:[-r,\infty)\mapsto
\mathbb{C}.$

$\widetilde{f}(z_t)$ is, due to (\ref{zt}), a function of
$u,\,\overline{u},$ and we develop it in a series of powers
\begin{equation}\label{seriesf}\widetilde{f}(z_t)=\sum_{j+k\geq2}\frac{1}{j!k!}f_{jk}(\alpha)u^{j}\overline{u}^{k}.
\end{equation}

Equation  (\ref{proj-eq1}) becomes
\begin{equation}\label{eq-series} \frac{du}{dt}=\widetilde{\lambda}_1(\alpha)
u+\sum_{j+k\geq2}\frac{1}{j!k!}g_{jk}(\alpha)u^{j}\overline{u}^{k},
\end{equation}
where
\begin{equation}\label{gjk}
g_{jk}(\alpha)=\Psi_1(0)f_{jk}(\alpha).\end{equation}
If the first Lyapunov coefficient at $\alpha^*$, \cite{K},
\begin{equation}\label{l1}
l_1(\alpha^*)=\frac{1}{2\omega^{*2}}Re(ig_{20}(\alpha^*)g_{11}(\alpha^*)+
\omega^* g_{21}(\alpha^*))
\end{equation}
is not equal to zero, then at $\alpha=\alpha^*$ our equation
presents a non-degenerated Hopf bifurcation.

We assume that the initial function, $\phi,$ is on the invariant
manifold. We give the algorithm for computing the first Lyapunov
coefficient. We are interested in the nonlinear part of our
equation. We consider the function $\beta(x)x,$ compute its
derivatives of order higher than one at the point $x_2,$ and denote
them  by $B_n$
\[B_n=\beta^{(n)}(x_2)x_2+n\beta^{(n-1)}(x_2). \\
\]
Let $F$ be the function $F(z)=\frac{1}{2!}B_2z^2+\frac{1}{3!}B_3
z^3+\frac{1}{4!}B_4z^4+...$. The nonlinear part of the RHS of
(\ref{nonlineq1}) is $f(z(t),z(t-r))=-F(z(t))+kF(z(t-r)),$ or, in
terms of $z_t,$ the function $\widetilde{f}$ of (\ref{nonlineq2}) is
$\widetilde{f}(z_t)=-F(z_t(0))+kF(z_t(-r)).$

The normal form of the restriction of the problem to the unstable
manifold is \cite{K}
\begin{equation}\label{nform}\frac{d}{dt}u=(\beta(\alpha)+i)u+su|u|^2,
\end{equation}
where
$\beta(\alpha)=\mu(\alpha)/\omega(\alpha),\,s=sign(l_1(\alpha^*)). $
For $\alpha=\alpha^*,\,\beta(\alpha^*)=0,$ and the normal form on
the center manifold is obtained. The existence and stability
properties of the cycle emerged by Hopf bifurcation are given by the
signs of $\beta(\alpha)$ (for $\alpha\in U^*$) and $l_1(\alpha^*).$
We must compute $l_1(\alpha^*),$ given by (\ref{l1}).

In order to identify the coefficients $g_{jk}(\alpha^*)$ we must
find the
 coefficients for the series of
$\widetilde{f}.$ The following computations are done at
$\alpha=\alpha^*,$ but for simplicity, the parameter will not be
written.

We have (by denoting $\varphi_1=\varphi$)
\begin{equation}\label{termnl}\widetilde{f}(z_t)=-F([\widetilde{T}(t)\phi](0))+kF([\widetilde{T}(t)\phi](-r))=
\end{equation}
\[=-\frac{1}{2!}B_2[u\varphi(0)+\overline{u}\overline{\varphi}(0)+\frac{1}{2}w_{20}(0)u^{2}
+
w_{11}(0)u\overline{u}+\frac{1}{2}w_{02}(0)\overline{u}^{2}+...]^2-
\]
\[-\frac{1}{3!}B_3[u\varphi(0)+\overline{u}\overline{\varphi}(0)
+\frac{1}{2}w_{20}(0)u^{2}+
w_{11}(0)u\overline{u}+\frac{1}{2}w_{02}(0)\overline{u}^{2}+...]^3-....+
\]
\[+\frac{1}{2!}kB_2[u\varphi(-r)+\overline{u}\overline{\varphi}(-r)+
\frac{1}{2}w_{20}(-r)u^{2}+
w_{11}(-r)u\overline{u}+\frac{1}{2}w_{02}(-r)\overline{u}^{2}+...]^2+
\]
\[+\frac{1}{3!}kB_3[u\varphi(-r)+\overline{u}\overline{\varphi}(-r)+
\frac{1}{2}w_{20}(-r)u^{2}+
w_{11}(-r)u\overline{u}+\frac{1}{2}w_{02}(-r)\overline{u}^{2}+...]^3+...=
\]
\[=\sum_{j+k\geq 2}\frac{1}{j!k!}f_{jk}u^j\overline{u}^k.
\]
The second order terms are
\[-\frac{1}{2}B_2(u^2+2u\overline{u}+\overline{u}^2)+
\frac{kB_2}{2}(u^2\varphi^2(-r)+2u\overline{u}+\overline{u}^2\overline{\varphi}^2(-r)),
\]
and the coefficients
\[f_{20}=B_2(k\varphi^2(-r)-1)=-B_2(1-ke^{-2\omega^* ir}),
\]
\[f_{11}=B_2(k-1),
\]
\[f_{02}=B_2(k\overline{\varphi}^2(-r)-1)=-B_2(1-ke^{2\omega^* ir}),
\]
follow. From these relations and \eqref{gjk}, we already find
$g_{20},\,g_{11},\,g_{02}.$

In order to determine $l_1$ we still need to compute the coefficient
$g_{21}.$ From (\ref{termnl})  we obtain
\begin{equation}\label{f21}
f_{21}=B_2\left(-2w_{11}(0)-w_{20}(0)+2ke^{-i\omega
r}w_{11}(-r)+ke^{i\omega r}w_{20}(-r)\right)-B_3\left(1-ke^{-i\omega
r}\right),
\end{equation}
hence we have to determine
$w_{20}(0),\,w_{20}(-r),\,w_{11}(0),\,w_{11}(-r).$

The functions $w_{jk}\in \mathcal{B}_C$ are determined from the
differential equations obtained by equating the same degree terms of
 (\cite{N}, \cite{I04}):
\begin{equation}\label{gendiffeq}\frac{\partial}{\partial s}\sum_{j+k\geq2}\frac{1}{j!k!}w_{jk}(s)u^{j}\overline{u}^{k}=
\sum_{j+k\geq2}\frac{1}{j!k!}g_{jk}u^{j}\overline{u}^{k}\varphi_1(s)+
\end{equation}
\[+
\sum_{j+k\geq2}\frac{1}{j!k!}\overline{g}_{jk}\overline{u}^{j}u^{k}\varphi_2(s)+\frac{\partial}{\partial
t}\sum_{j+k\geq2}\frac{1}{j!k!}w_{jk}(s)u^{j}\overline{u}^{k},
\]
and the integration constants are determined from \footnote{In the
first form of this article, published in Journal of Middle Volga
Mathematical Society, (journal formerly named Proceedings of ...), 
unhappily, relation \eqref{conddiffeq} was written in an erroneous
form. The subsequent computations were affected by the error
therein. The small corrections that we brought to our work in its
second version posted on arXiv did not cover this error. We make
here the necessary corrections, and we apologize to the readers that
might have been confused by our previous mistakes. An Erratum will
be send to Journal of Middle Volga Math. Soc. }
\begin{equation}\label{conddiffeq}
\frac{d}{dt}\sum_{j+k\geq2}\frac{1}{j!k!}w_{jk}(0)u^{j}\overline{u}^{k}+
\sum_{j+k\geq2}\frac{1}{j!k!}g_{jk}u^{j}\overline{u}^{k}\varphi_1(0)+
\sum_{j+k\geq2}\frac{1}{j!k!}\overline{g}_{jk}\overline{u}^{j}u^{k}\varphi_2(0)=
\end{equation}
\[-(B_1+\delta)\sum_{j+k\geq2}\frac{1}{j!k!}w_{jk}(0)u^{j}\overline{u}^{k}+kB_1\sum_{j+k\geq2}\frac{1}{j!k!}w_{jk}(-r)u^{j}\overline{u}^{k}+
\sum_{j+k\geq2}\frac{1}{j!k!}f_{jk}u^{j}\overline{u}^{k}
\]

Hence, the function $w_{20}(\cdot)$ is solution of the differential
equation,
\[w'_{20}(s)=2\omega^* i
w_{20}(s)+g_{20}\varphi(s)+\overline{g}_{02}\overline{\varphi}(s).\]
We integrate between $-r$ and $0$ and obtain
\[w_{20}(0)-w_{20}(-r)e^{2\omega^* ir}=
\frac{g_{20}i}{\omega^*}(1-e^{\omega^*
ir})+\overline{g}_{02}\frac{i}{3\omega^*}(1-e^{3\omega^* ir}).
\]
By equating the coefficients of  $u^2$ in \eqref{conddiffeq}, the
relation
\[2\omega^* i w_{20}(0)+g_{20}+\overline{g}_{02}=-(B_1+\delta)w_{20}(0)+kB_1w_{20}(-r)+f_{20}
\]
results. From the two relations above, we get
\[w_{20}(0)=c \left[e^{2\omega^* ir}f_{20}+g_{20}\left(-\frac{kB_1 i}{\omega^*}+\frac{kB_1 i}{\omega^*}
e^{i\omega^* r}-e^{2i\omega^* r} \right)+\right.\]
\[\left.+\overline{g}_{02}\left(-\frac{kB_1i}{3\omega^*}+\frac{kB_1i}{3\omega^*}e^{3i\omega^* r}-e^{2i\omega^*
r}\right)\right],
\]
\[w_{20}(-r)=c \left[f_{20}+g_{20}(1-2e^{\omega^*
ir}-\frac{B_1+\delta}{\omega^*}i(1-e^{\omega^{*} i r}))-\right.\]
\[\left.-\frac{1}{3}\overline{g}_{02}(1+2e^{3\omega^*
ir}+\frac{B_1+\delta}{\omega^*}i(1-e^{3\omega^{*} i r}))\right],
\]
where
$c=[-kB_1+(B_1+\delta)\cos(2\omega^*r)-2\omega^*\sin(2\omega^*r)-i(2\omega^*\cos(2\omega^*r)+(B_1+\delta)\sin(2\omega^*r))]/
[(kB_1)^2+(B_1+\delta)^2+(2\omega^*)^2-2kB_1(B_1+\delta)\cos(2\omega^*r)+4kB_1\omega^*\sin(2\omega^*r)].$

We now compute $w_{11}.$ It is the solution of
\[\frac{d}{ds}w_{11}(s)=g_{11}e^{\omega^* is}+\overline{g}_{11}e^{-\omega^*
is}.
\]

By integrating between $-r$ and 0, we obtain
\[w_{11}(0)=w_{11}(-r)-\frac{i}{\omega^* }g_{11}(1-e^{-\omega^* ir})+\frac{i}{\omega^* }\overline{g}_{11}(1-e^{\omega^* ir}).
\]

The second relation,  is here
\[g_{11}+\overline{g}_{11}=
-(B_1+\delta)w_{11}(0)+kB_1w_{11}(-r)+f_{11}.
\]
It follows that
\[w_{11}(0)=c_1\left[f_{11}-g_{11}-\overline{g}_{11}+\frac{kB_1i}{\omega^*}\left(g_{11}(1-e^{-\omega^* i r})
-\overline{g}_{11}(1-e^{\omega^*ir})\right)\right],
\]
\[w_{11}(-r)=c_1\left[f_{11}-g_{11}-\overline{g}_{11}+\frac{(B_1+\delta)i}{\omega^*}\left(g_{11}(1-e^{-\omega^* i r})
-\overline{g}_{11}(1-e^{\omega^*ir})\right)\right].
\]
with $c_1=1/(B_1+\delta-kB_1).$
 Now we have all the elements
for determining $f_{21}$ by formula (\ref{f21}), and then we can
compute
\[g_{21}=\Psi_1(0)f_{21},
\]
and after that, the first Lyapunov coefficient, that determines the
stability properties of the periodic orbit emerged by Hopf
bifurcation.

The above indicated computations, are easy to perform with a
computing program (Maple e.g.), for a given Hopf bifurcation point,
i.e. if all the parameters are known. If we try to make the
computations with all (or some of) the parameters non-determined, we
are lead to very complicated expressions that can not be handled
even by the program Maple itself.

\subsection{An example}
In order to provide an example, we choose the following "strategy"
for finding Hopf bifurcation points. We choose
$n^*,\,\beta^*_0,\,\delta^*, k^*,$ in an acceptable from biological
point of view zone of parameters. With these values of the
parameters we can compute $B^*_1$, (the value of $B_1$ at these
parameters values), $p^*=\delta^*+B^*_1,\,q^*=k^*B^*_1,\,$ and
determine $\omega^*$ and $r^*$ such that a Hopf bifurcation actually
takes place. That is, we set $\omega^*=\sqrt{(q^*)^2-(p^*)^2}$ and
$r^*=\arccos(p^*/q^*)/\omega^*.$ Once $r^*$ determined, we can also
compute $\gamma^*$.

As before, we choose one of the parameters $a_j$ to be the variable
one and check the condition that the curve $\lambda_1(a_j)$ actually
crosses the imaginary axis in the complex plane, when $a_j$ varies
around $a_j^*$. This can be done by example by differentiating
relations
 (\ref{re,im}) with respect to $a_j$ and by solving the obtained
 system of equations with respect to
$\mu'(a_j^*),\,\omega'(a_j^*)$.

\textbf{The example.}  We consider the following selection of
parameters:
\[n^*=12,\,\beta^*_0=1.77,\,\delta^*=0.05,\,k^*=1.180746972\, (\Leftrightarrow
\gamma^* r^*=0.527).
\]
These parameters were chosen in the zone of biological interest, as
can be seen in \cite{PM-M}.

We find then $x_2=1.150859618,$
$B^*_1=-2.524121872,$\,$p^*=\delta^*+B^*_1=-2.474121872,$\,$q^*=k^*B^*_1=-2.980349005,\,$
$\omega^*=\sqrt{q^{*2}-p^{*2}}=1.661686238,$\,\\$r^*=\arccos(p^*/q^*)/\omega^*=.3559207407,$\,
$\gamma^*=1.48067.$

We choose to vary $r$ and, by differentiating the relations
(\ref{re,im}) with respect to $r$ find that $\mu'(r^*)=25.66.$

The first Lyapunov coefficient is easily computed since all the
parameters are known, by following the algorithm in the previous
subsection. We get
\[1_1(r^*)=-43.71063.
\]
Thus the normal form of the restriction of the equation to the
unstable and to the center manifold is (by (\ref{nform}))
\[\frac{d}{dt}u=(\beta(r)+i)u-u|u|^2,
\]
with $\beta(r)=\frac{\mu(r)}{\omega(r)}$ and $\beta(r^*)=0.$ The
normal form in polar coordinates is
$$
\left\{%
\begin{array}{ll}
     & \dot{\rho}=(\beta(r)-\rho^2)\rho,
 \\
     & \dot{\theta}=1. \\
\end{array}%
\right.
$$
We are interested in the behavior of the solutions on a small
neighborhood of $r^*.$ Since $\mu'(r^*)>0,$ on such a neighborhood
 $r<r^*\,\Rightarrow \mu(r)<0\,\Rightarrow\,\beta(r)<0,$ and the
 equation in $\rho$ above has the only
 equilibrium solution $\rho=0\,\Leftrightarrow\,x=x_2$ and it is
 stable, while  $r>r^*\,\Rightarrow\, \mu(r)>0\,\Rightarrow\, \beta(r)>0.$
Thus a stable periodic solution exists for $r>r^*$.

\begin{center}

\includegraphics[width=7cm,height=3cm]{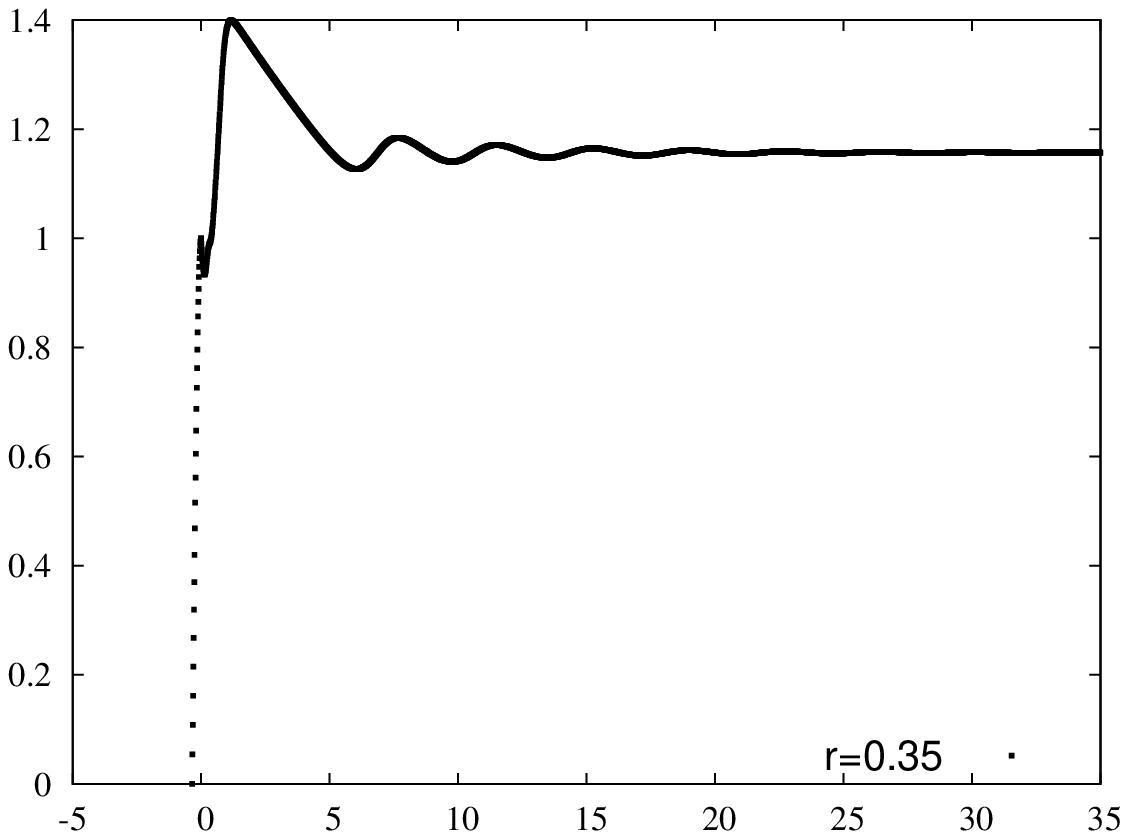}\newline
\textsf{a)  r=0.35.}

\end{center}

\begin{center}

\includegraphics[width=7cm,height=3cm]{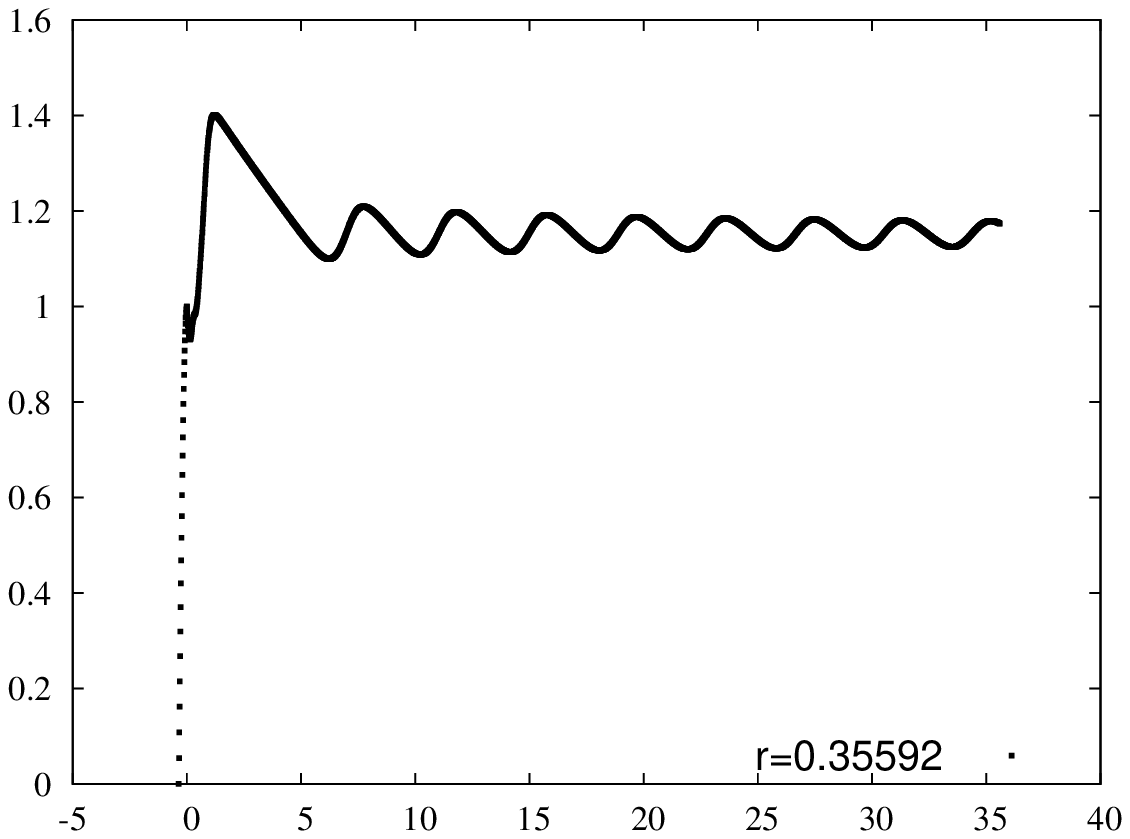}\newline
\textsf{b) r=0.35592.}

\end{center}

\begin{center}

\includegraphics[width=7cm,height=3cm]{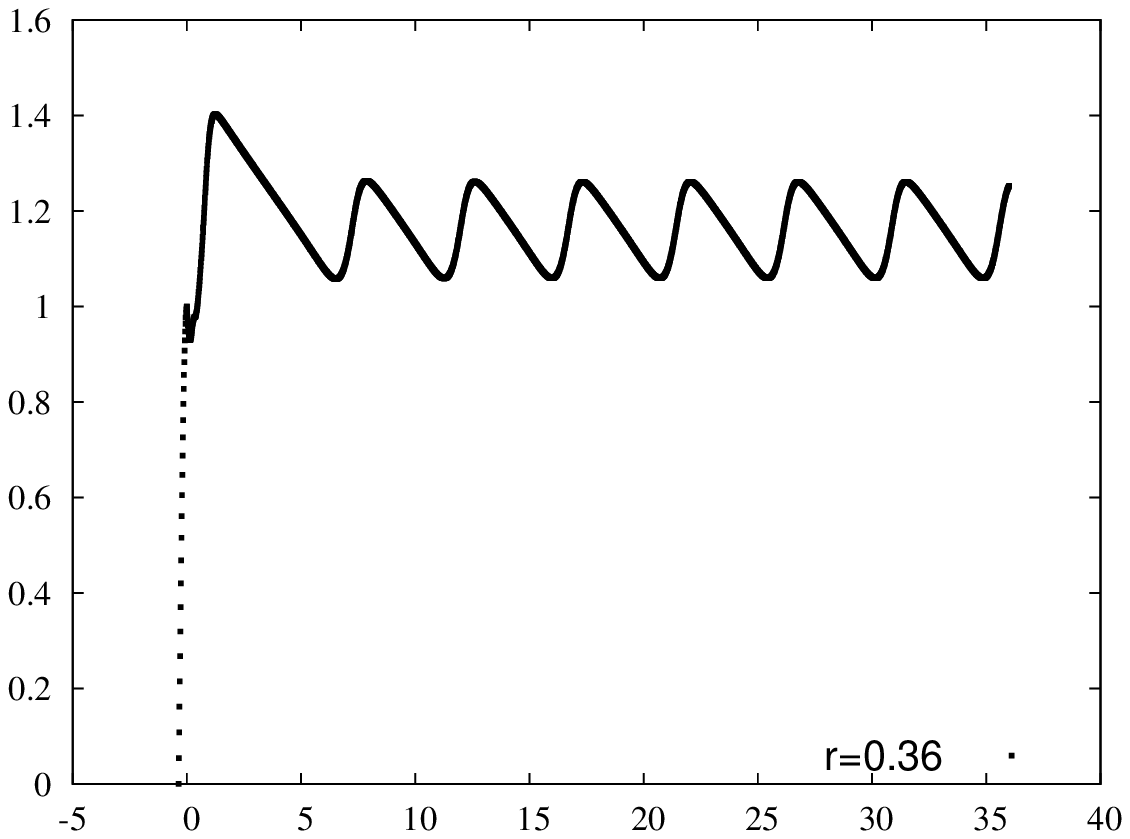}\newline
\textsf{c)  r=0.36.\\Fig.1) Trajectories obtained by numerical
integration of eq. (\ref{eq}) for initial condition
$\phi(s)=cos(\frac{\pi}{2r}s),$ for some values of $r$ around the
Hopf bifurcation value $r^*=.3559207407$. }
\end{center}

These assertions seem in contradiction with the condition of
stability (\ref{A}). Actually there is no contradiction.  The
function $g$ defined in (\ref{h}) satisfies, for the chosen
parameters values, $g(r^*)=0$ and $\frac{dg}{dr}(r^*)=-20.236.$
Thus, $r<r^*\,\Rightarrow\,g(r)>0$ and the equilibrium point $x_2$
is stable, while for $r>r^*$, $g(r)<0$ and the equilibrium point is
unstable.

The numerical integration confirms our theoretical analysis, as the
figures above show.

\end{document}